\documentstyle{article}
 
\normalbaselines 
\input amssym.def
\input amssym.tex

\newcount\skewfactor
\def\mathunderaccent#1#2 {\let\theaccent#1\skewfactor#2
\mathpalette\putaccentunder}
\def\putaccentunder#1#2{\oalign{$#1#2$\crcr\hidewidth
\vbox to.2ex{\hbox{$#1\skew\skewfactor\theaccent{}$}\vss}\hidewidth}}
\def\name{\mathunderaccent\tilde-3 }

\newcommand{\forces}{\Vdash}

\newcommand{\V}{{\bf V}} 
\newcommand{\lesdot}{\mathrel{\mathord{<}\!\!\raise 0.8
pt\hbox{$\scriptstyle\circ$}}}  
\newcommand{\conc}{{}^\frown\!}
\newcommand{\lh}{\ell g\/} 
\newcommand{\rest}{{\restriction}}
\newcommand{\dom}{{\rm dom}} 

\newcommand{\cf}{{\rm cf}}
\newcommand{\bP}{{\Bbb P}}
\newcommand{\bR}{{\Bbb R}}
\newcommand{\bQ}{{\Bbb Q}}
\newcommand{\bfu}{{\bf u}}
\newcommand{\nbQ}{{\name{\Bbb Q}}}
\newcommand{\QED}{\hfill\vrule width 6pt height 6pt depth 0pt 
\vspace{0.1in}} 
\newcommand{\Proof}{\noindent{\sc Proof} \hspace{0.2in}} 

\newcommand{\D}{{\cal D}}
\newcommand{\cQ}{{\cal Q}}
\newcommand{\nD}{\name{\D}}
\newcommand{\lx}{{\ell x}}
\newcommand{\can}{{}^{\textstyle \omega}2} 
\newcommand{\fs}{{}^{\textstyle <\!\omega}2} 
\newcommand{\baire}{{}^{\textstyle \omega}\omega} 
\newcommand{\iso}{[\omega]^{\textstyle \aleph_0}} 
 
\newcommand{\fseo}{{}^{\textstyle <\!\omega}\omega} 

\newcommand{\qdwa}{\bQ^2_{I,h}}
\newcommand{\qtrzy}{\bQ^3_{I,h}}
\newcommand{\qcztery}{\bQ^4_{I,h}}
\newcommand{\qjeden}{\bQ^1_{I,h}}
\newcommand{\qell}{\bQ^\ell_{I,h}}
\newcommand{\qnic}{\bQ_{I,h}}

\newcommand{\var}{{\rm var}}
\newcommand{\nwd}{{\rm nwd}}

\newtheorem{theorem}{Theorem}[section] 

\newtheorem{proposition}[theorem]{Proposition} 
\newtheorem{corollary}[theorem]{Corollary} 
\newtheorem{conclusion}[theorem]{Conclusion} 
\newtheorem{definition}[theorem]{Definition}
\newtheorem{remark}[theorem]{Remark}
\newtheorem{lemma}[theorem]{Lemma}
\newtheorem{concrem}[theorem]{Concluding Remarks}
\title{There may be no nowhere dense ultrafilter}
\author{{\bf Saharon Shelah}\thanks{\ The research partially supported by
``Basic Research Foundation'' of the Israel Academy of Sciences and
Humanities. Publication 594.}\\
Institute of Mathematics\\
The Hebrew University of Jerusalem\\
91904 Jerusalem, Israel\\
and\\
Department of Mathematics\\
Rutgers University\\
New Brunswick, NJ 08854, USA}
\date{\today}

\begin{document} 
\maketitle 

\begin{abstract}
We show the consistency of ZFC + "there is no NWD-ultrafilter on
$\omega$", which means: for every non-principal ultrafilter $\D$ on the set of
natural numbers, there is a function $f$ from the set of natural numbers to
the reals, such that for every nowhere dense set $A$ of reals, $\{n: f(n)\in
A\}\notin\D$. This answers a question of van Douwen, which was put in more
general context by Baumgartner.
\end{abstract}
\vfill
\eject

\section{Introduction}
We prove here the consistency of ``there is no NWD-ultrafilter on $\omega$"
(non-principal, of course). This answers a question of van Douwen \cite{vD81}
which appears as question 31 of \cite{B6}. Baumgartner \cite{B6} considers the
question which he dealt more generally with $J$-ultrafilter where
\begin{definition}  
\label{0.1}
\begin{enumerate}
\item An ultrafilter $\D$, say on $\omega$, is called a $J$-ultrafilter where
$J$ is an ideal on some set $X$ (to which all singletons belong, to avoid
trivialities) {\em if} for every function $f:\omega\longrightarrow X$ for some
$A \in \D$ we have $f''(A)\in J$. 
\item  The NWD--ultrafilters are the $J$-ultrafilters for
$J=\{B\subseteq\cQ:B$ is nowhere dense$\}$ ($\cQ$ is the set of all rationals;
we will use an equivalent version, see \ref{2.2}).
\end{enumerate}
\end{definition}
This is also relevant for the consistency of ``every (non-trivial) c.c.c.
$\sigma$-centered forcing notion adds a Cohen real", see \cite{Sh:F151}.

The most natural approach to a proof of the consistency of ``there is no
NWD--ultrafilter'' was to generalize the proof of CON(there is no $P$-point)
(see \cite[VI, \S4]{Sh:b} or \cite[VI, \S4]{Sh:f}), but I (and probably
others) have not seen how. 

We use an idea taken from \cite{Sh:407}, which is to replace the given maximal
ideal $I$ on $\omega$ by a quotient; moreover, we allow ourselves to change 
the quotient. In fact, the forcing here is simpler than the one in
\cite{Sh:407}. A related work is Goldstern Shelah \cite{GoSh:388}. 

We similarly may consider the consistency of ``no $\alpha$--ultrafilter" for
limit $\alpha < \omega_1$ (see \cite{B6} for definition and discussion of
$\alpha$--ultrafilters). This question and the problems of preservation of 
ultrafilters and distinguishing existence properties of ultrafilters will be
dealt with in a subsequent work \cite{Sh:F187}.  

In \S3 we note that any ultrafilter with property $M$ (see Definition
\ref{2.5}) is an NWD--ultrafilter, hence it is consistent that there is no
ultrafilter (on $\omega$) with property $M$. 

I would like to thank James Baumgartner for arousing my interest in
the questions on NWD--ultrafilters and $\alpha$-ultrafilters and
Benedikt on asking about the property $M$ as well as Shmuel Lifches
for corrections, the participants of my seminar in logic in Madison
Spring'96 for hearing it, and Andrzej Ros{\l}anowski for corrections and
introducing the improvements from the lecture to the paper. 

\section{The basic forcing} 
In Definition \ref{1.1} below we define the forcing notion $\qjeden$ which
will be the one used in the proof of the main result \ref{2.4}. The other
forcing notion defined below, $\qdwa$, is a relative of $\qjeden$. Various
properties are much easier to check for $\qdwa$, but unfortunately it does not
do the job. The reader interested in the main result of the paper only may 
concentrate on $\qjeden$.

\begin{definition}
\label{Iequiv}
Let $I$ be an ideal on $\omega$ containing the family $[\omega]^{\textstyle
{<}\omega}$ of finite subsets of $\omega$.
\begin{enumerate}
\item We say that an equivalence relation $E$ is an {\em $I$--equivalence
relation} if: 
\begin{description}
\item[(a)]  $\dom(E)\subseteq \omega$,
\item[(b)]  $\omega \setminus \dom(E) \in I$,
\item[(c)]  each $E$-equivalence class is in $I$.
\end{description}
\item For $I$-equivalence relations $E_1,E_2$ we write $E_1 \leq E_2$ {\em if}
\begin{description}
\item[(i)]\ \ $\dom(E_2) \subseteq \dom(E_1)$,
\item[(ii)]\  $E_1\restriction\dom(E_2)$ refines $E_2$,
\item[(iii)]  $\dom(E_2)$ is the union of a family of $E_1$-equivalence
classes. 
\end{description}
\end{enumerate}
\end{definition}

\begin{definition}
\label{1.1}
Let $I$ be an ideal on $\omega$ to which all finite subsets of $\omega$ belong
and let $h:\omega\longrightarrow\omega$ be a non-decreasing function. Let
$\ell\in\{1,2\}$. We define a forcing notion $\qell$ (if $h(n) = n$ we may 
omit it) intended to add $\langle y^n_i:i < h(n),n <\omega\rangle$, $y^n_i\in
\{-1,1\}$. We use $x^n_i$ as variables.
\begin{enumerate}
\item $p\in\qell$ if and only if $p=(H,E,A)=(H^p,E^p,A^p)$ and
\begin{description}
\item[(a)]  $E$ is an $I$--equivalence relation on $\dom(E)\subseteq \omega$,
\item[(b)]  $A = \{n\in\dom(E): n =\min(n/E)\}$,
\item[(c)]  if $\ell = 1$, {\em then} $H$ is a function with range $\subseteq
\{-1,1\}$ and domain 
\[\begin{array}{ll}
B^p_1=\big\{x^n_i:& i<h(n)\mbox{ and } n\in\omega\setminus\dom(E)\ \mbox{\em
  or}\\  
\ &n\in\dom(E)\mbox{ and }i\in [h(\min(n/E)),h(n))\big\},
  \end{array}\]
\item[(d)] if $\ell=2$, {\em then} 
\begin{description}
\item[$(\alpha)$]  $H$ is a function on $\dom(H)=B^p_2\cup B^p_3$, where
\[\hspace{-1cm}\begin{array}{lr}
B^p_2=\{x^m_i:m\in\omega,\ A^p\cap (m+1)=\emptyset,\ i<h(m)\} &\mbox{ and
}\\
B^p_3 = \{x^m_i:m\in\dom(E^p)\setminus A^p\mbox{ or } m\notin\dom(E^p)
\mbox{ but }&A^p\cap m\neq\emptyset,\ \ \\
&i < h(m)\},
  \end{array}\]
\item[$(\beta)$]  for $x^m_i \in B^p_3$, $H(x^m_i)$ is a function of the
variables $\{x^n_j:(n,j) \in w_p(m,i)\}$ to $\{-1,1\}$, where 
\[w_p(m) = w_p(m,i) = \{(\ell,j):\ell\in A^p\cap m \mbox{ and } j<h(\ell)\},\]
for $n \in A^p$ we stipulate $H^p(x^n_i) = x^n_i$ and
\item[$(\gamma)$] $H\restriction B^p_2$ is a function to $\{-1,1\}$.
\end{description}
\item[(e)] if $\ell = 2$ and $x^n_i\in B^p_3$, $n^*=\min(n/E^p)<n$ and $y^m_i
\in\{-1,1\}$ for $m\in A^p \cap n^*$, $i<h(m)$ and $z^n_j\in\{-1,1\}$ for $j< 
h(n^*)$ {\em then} for some $y^{n^*}_j\in \{-1,1\}$ for $j<h(n^*)$ we have
\[j< h(n^*)\quad\Rightarrow\quad z^n_j=(H^p(x^n_j))(\ldots,y^m_i,\ldots)_{(m,
i) \in w_p(n,j)}.\]
\end{description}
When it can not cause any confusion, or we mean ``for both $\ell=1$ and
$\ell=2$'', we omit the superscript $\ell$.  
\item  Defining functions like $H(x^m_i)$, $x^m_i\in B^p_3$ (when $\ell = 2$), 
we may allow to use dummy variables. In particular, if $H^p(x^m_i)$ is $-1,1$
we identify it with constant functions with this value.  
\item We say that a function $f:\{x^n_i:i< h(n),\ n<\omega\}\longrightarrow 
\{-1,1\}$ satisfies a condition $p\in\qell$ if:
\begin{description}
\item[(a)]  $f(x^n_i) = H^p(x^n_i)$ {\em when} $x^n_i\in B^p_1$ and $\ell=1$,
or $x^n_i\in B^p_2$ and $\ell=2$, 
\item[(b)]  $f(x^n_i) = H^p(x^n_i)(\ldots,f(x^m_j),\ldots)_{(m,j)\in w_p(n,
i)}$ {\em when} $\ell=2$ and $x^n_i\in B^p_3$,
\item[(c)]  $f(x^n_i) = f(x^{\min(n/E^p)}_i)$ {\em when} $\ell = 1$, $n\in
\dom(E^p)$ and $i<h(\min(n/E^p))$.
\end{description}
\item The partial order $\leq=\leq_{\qell}$ is defined by $p\leq q$ if and
only if:
\begin{description}
\item[$(\alpha)$]  $E^p\leq E^q$,
\item[$(\beta)$]   every function $f:\{x^n_i:i<h(n),\ n<\omega)\}
\longrightarrow\{-1,1\}$ satisfying $q$ satisfies $p$.
\end{description}
\end{enumerate}
\end{definition}

\begin{proposition}
\label{1.2} 
$(\qell,\leq_{\qell})$ is a partial order. \QED
\end{proposition}

\begin{remark}
{\em
We may reformulate the definition of the partial orders $\leq_{\qell}$, making
them perhaps more direct. Thus, in particular, if $p,q\in\qjeden$ then
$p\leq_{\qjeden} q$ if and only if the demand $(\alpha)$ of \ref{1.1}(4) holds
and 
\begin{description}
\item[$(\beta)^*$] for each $x^n_i\in B^q_1$:
\begin{description}
\item[(i)]  \ \ if $x^n_i\in B^p_1$ then $H^q(x^n_i)=H^p(x^n_i)$,
\item[(ii)] \ if $n\in\dom(E^p)\setminus\dom(E^q)$, $i<h(\min(n/E^p))$ then
$H^q(x^n_i)=H^q(x^{\min(n/E^p)}_i)$, 
\item[(iii)]  if $n\in\dom(E^q)\setminus A^p$, $\min(n/E^p)>\min(n/E^q)$ and
$h(\min(n/E^q))\leq i<h(\min(n/E^p))$ then
$H^q(x^n_i)=H^q(x^{\min(n/E^p)}_i)$. 
\end{description}
\end{description}
The corresponding reformulation for the forcing notion $\qdwa$ is more
complicated, but it should be clear too.
}
\end{remark}

One may wonder why we have $h$ in the definition of $\qell$ and we do not fix
that e.g. $h(n)=n$. This is to be able to describe nicely what is the forcing
notion $\qell$ below a condition $p$ like. The point is that $\qell\rest\{q:
q\geq p\}$ is like $\qell$ but we replace $I$ by its quotient and we change
the function $h$. More precisely:

\begin{proposition}
\label{1.3}
If $p\in\qell$ and $A^p=\{n_k:k<\omega\}$, $n_k<n_{k+1}$, $h^*:\omega
\longrightarrow \omega$ is $h^*(k)=h(n_k)$ and $I^*=\{B\subseteq\omega:
\bigcup\limits_{k\in B}(n_k/E) \in I\}$ {\em then} $\qell\restriction
\{q:p \leq_{\qell} q\}$ is isomorphic to $\bQ^\ell_{I^*,h^*}$.
\end{proposition}

\Proof  Natural. \QED

\begin{definition}
\label{1.4}
We define a $\qnic$--name $\name{\bar{\eta}}=\langle\name{\eta}_n:n<\omega
\rangle$ by:\\
$\name{\eta}_n$ is a sequence of length $h(n)$ of members of $\{-1,1\}$ such
that 
\[\name{\eta}_n[G_{\qnic}](i)=1\quad\Leftrightarrow\quad (\exists p\in
G_{\qnic})(H^p(x^n_i)=1).\]
[Note that in both cases $\ell=1$ and $\ell=2$, if $H^p(x^n_i)=1$, $x^n_i\in
\dom(H^p)$ and $q\geq p$ then $H^q(x^n_i)=1$; remember \ref{1.1}(2).]
\end{definition}

\begin{proposition}
\label{1.5}
\begin{enumerate}
\item If $n<\omega$, $A^p\cap(n+1)=\emptyset$ {\em then} $p\forces\mbox{``}
\name{\eta}_n=\langle H^p(x^n_i):i<h(n)\rangle$''.
\item For each $n<\omega$ the set $\{p\in\qnic: A^p \cap(n+1)=\emptyset\}$ is
dense in $\qnic$. 
\item If $p\in {\Bbb Q}_{I,h}$ and $a\subseteq A^p$ is finite or at least
$\bigcup\limits_{n\in a}(n/E^p)\in I$, and 
\[f:\{x^n_i:i < h(n) \mbox{ and }n \in a\} \longrightarrow \{-1,1\},\]
{\em then} for some unique $q$ which we denote by $p^{[f]}$, we have:
\begin{description}
\item[(a)]  $p \leq q \in \qnic$,
\item[(b)]  $E^q = E^p \restriction \bigcup\{n/E^p:n \in A \setminus a\}$,
\item[(c)]  for $n\in a$, $i<h(n)$ we have $H^q(x^n_i)$ is $f(x^n_i)$.
\end{description}
\end{enumerate}
\end{proposition}

\Proof Straight. \QED

\begin{definition}
\label{1.6} 
\begin{enumerate}
\item  $p \leq_n q$ (in $\qnic$) {\em if} $p \le q$ and:
\[k \in A^p \ \&\  |A^p \cap k| < n\quad \Rightarrow\quad k\in A^q.\]
\item $p \le^*_n q$ {\em if} $p \le q$ and:
\[k \in A^p \ \&\  |A^p \cap k| < n \quad\Rightarrow\quad k\in A^q \ \&\ 
k/E^p = k/E^q.\]
\item $p\le^\otimes_n q$ {\em if} $p \le_{n+1} q$ and: 
\[n > 0\quad\Rightarrow\quad p \le^*_n q\qquad\mbox{ and }\quad\dom(E^q)=
\dom(E^p).\]
\item For a finite set $\bfu\subseteq\omega$ we let $\var(\bfu)\stackrel{\rm
def}{=}\{x^n_i:i < h(n),\ n \in\bfu\}$.
\end{enumerate}
\end{definition}

\begin{proposition}
\label{1.7}
\begin{enumerate}
\item If $p \le q$, $\bfu$ is a finite initial segment of $A^p$ and $A^q\cap
\bfu=\emptyset$, {\em then} for some unique $f:\{x^n_i:i< h(n)\mbox{ and } n
\in \bfu\}\longrightarrow \{-1,1\}$ we have $p \leq p^{[f]} \le q$ (where
$p^{[f]}$ is from \ref{1.5}(3)). 
\item If $p\in\qell$ and $\bfu$ is a finite initial segment of $A^p$ {\em
then} 
\begin{description}
\item[$(*)_1$]  $f \in {}^{\var(\bfu)}\{-1,1\}$ implies $p \le p^{[f]}$ and
\[p^{[f]} \Vdash\mbox{`` }(\forall n\in\bfu)(\forall i<h(n))(\name{\eta}_n(i)= 
f(x^n_i))\mbox{''},\] 
\item[$(*)_2$]  the set $\{p^{[f]}:f\in {}^{\var(\bfu)}\{-1,1\}\}$ is predense
above $p$ (in $\qell$). 
\end{description}
\item $\leq_n$ is a partial order on $\qell$, and $p \leq_{n+1} q\ \Rightarrow
\ p\leq_n q$.  Similarly for $<^*_n$ and $<^\otimes_n$.  Also
\[p\leq^\otimes_n q\quad \Rightarrow\quad p \leq^*_n q\quad \Rightarrow\quad p
\leq_n q\quad \Rightarrow\quad p \leq q.\]
\item If $p \in \qell$, $\bfu$ is a finite initial segment of $A^p$, $|\bfu|=
n$ and 
\[f:\{x^n_i:i < h(n) \mbox{ and } n \in \bfu\} \longrightarrow \{-1,1\}\quad
\mbox{ and }\quad p^{[f]} \leq q \in \qell,\]
{\em then} for some $r\in\qell$ we have $p\leq^*_n r \leq q$, $r^{[f]}=q$.
\item If $p\in \qdwa$, $\bfu$ is a finite initial segment of
$A^p$, $|\bfu| = n+1$ and 
\[f:\{x^n_i:i < h(n) \mbox{ and } n \in \bfu\} \longrightarrow \{-1,1\}\quad
\mbox{ and }\quad p^{[f]} \leq q,\]
{\em then} for some $r\in\qdwa$ we have $p<^\otimes_n r\leq q$ and $r^{[f]}=
q$. 
\end{enumerate}
\end{proposition}

\Proof  1)\ \ \ Define $f:\{x^n_i:i < h(n) \mbox{ and } n \in \bfu\}
\longrightarrow \{-1,1\}$ by: 
\[f(x^n_i)\mbox{ is the constant value of }H^q(x^n_i)\]
(it is a constant function by \ref{1.1}(1)(e), \ref{1.1}(1)(f($\gamma$))).

\noindent 2)\ \ \ By \ref{1.5} and \ref{1.7}(1).
  
\noindent 3)\ \ \ Check.

\noindent 4)\ \ \ First let us define the required condition $r$ in the case
$\ell=1$. So we let 
\[\begin{array}{l}
\dom(E^r) = \bigcup\limits_{n\in\bfu} (n/E^p)\cup \dom(E^q),\\
E^r=\big\{(n_1,n_2):\,n_1\,E^q\,n_2 \mbox{ or for some } n\in\bfu 
\mbox{ we have: }\{n_1,n_2\} \subseteq (n/E^p))\big\},\\
A^r = \bfu \cup A^q
  \end{array}\]
(note that if $n_1\,E^q\,n_2$ then $n_1 \notin \bfu$). Next, for $x^n_i\in
B^r_1$ (where $B^r_1$ is given by \ref{1.1}(1)(e)) we define 
\[H^r(x^n_i)=\left\{\begin{array}{ll}
H^q(x^n_i)& \mbox{ if } n\notin\bigcup\limits_{k\in\bfu}k/E^p\mbox{ and }
x^n_i \in\dom(H^q),\\
H^p(x^n_i)&\mbox{ if } n\in\bigcup\limits_{k\in\bfu} k/E^p\mbox{ and }
x^n_i\in \dom(H^p).
		    \end{array}\right.\]
It should be clear that $r=(H^r,E^r,A^r)\in\qjeden$ is as required.

If $\ell=2$ then we define $r$ in a similar manner, but we have to be more
careful defining the function $H^r$. Thus $E^r$ and $A^r$ are defined as
above, $B^r_2$, $B^r_3$ and $w_r(m,i)$ for $x^m_i\in B^r_3$ are given by
\ref{1.1}(1)(f). Note that $B^r_2=B^p_2$ and $B^r_3\subseteq B^p_3$. Next we
define:

if $x^m_i\in B^r_2$ then $H^r(x^m_i)=H^p(x^m_i)$,

if $x^m_i\in B^r_3$, $m\cap A^r\subseteq\bfu$ then $H^r(x^m_i)=H^p(x^m_i)$,

if $x^m_i\in B^r_3$ and $\min(\dom(E^q))<m$ then
\[\begin{array}{l}
H^r(x^m_i)(\ldots,x^k_j,\ldots)_{(k,j)\in w_r(m,i)}=\\
H^p(x^m_i)(x^k_j,H^q(x^{k'}_{j'})(\ldots,x^{k''}_{j''},\ldots)_{(k'',j'')\in
w_q(k',j')}))_{\scriptstyle (k,j)\in w_r(m,i)\atop \scriptstyle (k',j')\in
w_p(m,i)\setminus w_r(m,i)}.
  \end{array}\]
Note that if $(k',j')\in w_p(m,i)\setminus w_r(m,i)$, $x^m_i\in B^r_3$ then
$k'\in A^p\setminus (\bfu\cup A^q)$ and $w_q(k',j')\subseteq w_r(m,i)$. 

\noindent 5)\ \ \ Like the proof of (4). Let $n^*=\max(\bfu)$. Put $\dom(E^r)
=\dom(E^p)$ and declare that $n_1\, E^r\,n_2$ if one of the following occurs:
\begin{description}
\item[(a)] for some $n\in\bfu\setminus\{n^*\}$ we have $\{n_1,n_2\}\subseteq
(n/E^p)$, {\em or} 
\item[(b)] $n_1\,E^q\,n_2$ (so $n\in\bfu\Rightarrow \neg n\,E^p\,n_1$), {\em
or} 
\item[(c)] $\{n_1,n_2\}\subseteq B$, where
\[B\stackrel{\rm def}{=}n^*/E^p \cup \bigcup\{m/E^p:m\in\dom(E^p)\setminus
\dom(E^q),\ \min(m/E^p)>n^*\}.\]
\end{description}
We let $A^r = \bfu \cup A^q$ (in fact $A^r$ is defined from $E^r$). Finally
the function $H^r$ is defined exactly in the same manner as in (4) above (for
$\ell=2$). \QED

\begin{corollary}
\label{1.8}
If $p \in \qell$, $n<\omega$ and $\name{\tau}$ is a $\qell$--name of an
ordinal, {\em then} there are $\bfu,q$ and $\bar{\alpha}=\langle\alpha_f:f\in
{}^{\var(\bfu)}\{-1,1\} \rangle$ such that: 
\begin{description}
\item[(a)] $p \leq^*_n q \in \qell$,
\item[(b)] $\bfu = \{\ell \in A^p:| \ell \cap A^p| < n\}$,
\item[(c)] for $f \in {}^{\var(\bfu)}\{-1,1\}$ we have $q^{[f]}
\Vdash$ ``$\name{\tau} = \alpha_f"$,
\item[(d)] $q \Vdash$ ``$\name{\tau}\in\{\alpha_f:f\in {}^{\var(\bfu)}\{-1,
1\}\}"$ (which is a finite set).
\end{description}
\end{corollary}

\Proof Let $k =\prod\limits_{\ell\in\bfu}2^{h(\ell)}$. Let $\{f_\ell:\ell<k\}$
enumerate ${}^{\var(\bfu)}\{-1,1\}$. By induction on $\ell \le k$ define 
$r_\ell,\alpha_{f_\ell}$ such that:
\[r_0 = p,\quad r_\ell \leq^*_n r_{\ell + 1} \in \qell,\quad r^{[f_\ell]}_{
\ell+1}\Vdash_{\qell} \mbox{``}\name{\tau} = \alpha_{f_\ell}\mbox{''}.\]
The induction step is by \ref{1.7}(4). Now $q = r_k$ and $\langle \alpha_f:f
\in {}^{\var(\bfu)}\{-1,1\}\rangle$ are as required. \QED

\begin{corollary}
\label{1.8A}
If $\ell = 2$ then in \ref{1.8}(a) we may require $p\leq^\otimes_n
q\in\qell$.
\end{corollary}

\Proof  Similar: just use \ref{1.7}(5) instead of \ref{1.7}(4). \QED

\begin{definition}
\label{1.9}
Let $I$ be an ideal on $\omega$ containing $[\omega]^{\textstyle{<}\omega}$ and
let $E$ be an $I$--equivalence relation. 
\begin{enumerate}
\item We define a game $GM_I(E)$ between two players. The game lasts $\omega$
moves. In the $n^{\rm th}$ move the first player chooses an $I$-equivalence
relation $E^1_n$ such that 
\[E^1_0= E,\qquad [n > 0\quad \Rightarrow\quad E^2_{n-1} \leq E^1_n],\]
and the second player chooses an $I$-equivalence relation $E^2_n$ such that
$E^1_n \leq E^2_n$. In the end, the second player wins if
\[\bigcup \{\dom(E^2_n) \setminus\dom(E^1_{n+1}):n\in\omega\} \in I\]
(otherwise the first player wins).
\item For a countable elementary submodel $N$ of $({\cal
H}(\chi),{\in},{<}^*)$ such that $I,E\in N$ we define a game $GM^N_I(E)$ in a
similar manner as $GM_I(E)$, but we demand additionally that the relations
played by both players are from $N$ (i.e.~$E^1_n,E^2_n\in N$ for
$n\in\omega$). 
\end{enumerate}
\end{definition}

\begin{proposition}
\label{1.10} 
\begin{enumerate}
\item Assume that $I$ is a maximal (non-principal) ideal on $\omega$ and $E$
is an $I$--equivalence relation. Then the game $GM_I(E)$ is not
determined. Moreover, for each countable $N\prec({\cal H}(\chi),{\in},{<}^*)$
such that $I,E\in N$ the game $GM^N_I(E)$ is not determined.
\item For the conclusion of (1) it is enough to assume that ${\cal P}(\omega)/I
\models$ ccc.
\end{enumerate}
\end{proposition}

\Proof  1)\ \ \ As each player can imitate the other's strategy. 

\noindent 2)\ \ \ Easy, too, and will not be used in this paper. \QED

\begin{proposition}
\label{1.11}
\begin{enumerate}
\item Let $p\in\qell$. Suppose that the first player has no winning strategy in
$GM_I(E^p)$. {\em Then} in the following game Player I has no winning
strategy:  
\begin{description}
\item[] in the $n^{\rm th}$ move, 

Player I chooses a $\qell$-name $\name{\tau}_n$ of an ordinal and 

Player II chooses $p_n,\bfu_n,w_n$ such that: $w_n$ is a set of $\leq 
\prod\limits_{\ell\in \bfu_n} 2^{h(\ell)}$ ordinals, $p\leq p_n \leq^*_n
p_{n+1}$, $p_n\leq_{n+1} p_{n+1}$, $\bfu_n$ a finite initial segment of
$A^{p_n}$ with $n$ elements and $p_n \Vdash ``\name{\tau}_n\in w_n"$, moreover
\[f \in {}^{\var(\bfu_n)}\{-1,1\}\quad \Rightarrow\quad p^{[f]}_n\mbox{ forces
a value to }\name{\tau}_n.\]  
\item[] In the end, the second player wins if for some $q \ge p$ we have 
\[q \Vdash {`` }(\forall n\in\omega)(\name{\tau}_n \in w_n)\mbox{ ''.}\]
\end{description}
We can let Player II choose $k_n<\omega$ and demand $|\bfu_n| \le k_n$, and in
the end Player II wins if $\lim\inf\langle k_n:n < \omega \rangle < \omega$ or
there is $q$ as above. 
\item Let $p\in\qell$ and let $N$ be a countable elementary submodel of
$({\cal H}(\chi),{\in},{<}^*)$ such that $p,I,h\in N$. If the first player has
no winning strategy in $GM^N_I(E^p)$ then Player I has no winning strategy in
the game like above but with restriction that $\name{\tau}_n, p_n\in N$.
\end{enumerate}
\end{proposition}

\Proof 1)\ \ \ As in \cite[1.11, p.436]{Sh:407}.\\
Let ${\bf St}_p$ be a strategy for Player I in the game from \ref{1.11}. We
shall define a strategy ${\bf St}$ for the first player in $GM_I(E^p)$ during
which the first player, on a side, plays a play of the game from \ref{1.11}, using
${\bf St}_p$, with $\langle p_\ell: \ell< \omega \rangle$ and he also chooses
$\langle q_\ell:\ell<\omega\rangle$. 

Then, as ${\bf St}$ cannot be a winning strategy in $GM_I(E)$, in some play in
which the first player uses his strategy ${\bf St}$ he loses, and then
$\langle p_\ell:\ell < \omega \rangle$ will have an upper bound as required.

In the $n^{\rm th}$ move (so $E^1_\ell,E^2_\ell,q_\ell,p_\ell,\bfu_\ell,w_\ell$ 
for $\ell<n$ are defined), the first player in addition to choosing $E^1_n$
chooses $q_n,p_n,\bfu_n$, such that:
\begin{description}
\item[(a)]  $p=p_{-1}\leq q_0 = p_0$, $p_n \in \qell$, $q_n\in\qell$,
\item[(b)]  $p_n \le^*_n p_{n+1} \in \qell$,
\item[(c)]  $\bfu_0$ is $\emptyset$,
\item[(d)]  $\bfu_{n+1} = \bfu_n \cup\{\min(A^{q_{n+1}}\setminus \bfu_n)\}$,
so $|\bfu_{n+1}| = n+1$,
\item[(e)]  $E^1_0=E^p$, $E^1_{n+1}=E^{p_n} \restriction \big(\dom(E^{p_n})
\setminus \bigcup\limits_{i\in\bfu_n} i/E^{p_n}\big)$, 
\item[(f)]  $q_n$ is defined as follows:
\begin{description}
\item[$(f_0)$]  if $n=0$ then $E^{q_n}=E^2_0$,
\item[$(f_1)$]  if $n>0$ then $\dom(E^{q_n})=\dom(E^{p_{n-1}})$ and $x\,
E^{q_n}\, y$ if and only if

either $x\, E^2_n\, y$,

or for some $k\in\bfu_{n-1}$ we have $x,y\in k/E^{p_{n-1}}$'

or $x,y\in\big(\dom(E^1_n)\setminus\dom(E^2_n)\big)\cup\min(\dom(E^2_n))/
E^2_n$, 
\item[$(f_2)$]  $H^{q_n}$ is such that $p_{n-1}\leq q_n$,
\end{description}
\item[(g)] $p_n \leq^*_n q_{n+1} \leq^*_{n+1} p_{n+1}$, $p_n\leq_{n+1}
q_{n+1}$ (so $p_n\leq_{n+1} p_{n+1}$), 
\item[(h)] if $f \in {}^{\var(\bfu_n)}\{-1,1\}$ then $p^{[f]}_n$ forces a
value to $\name{\tau}_n$.
\end{description}
In the first move, when $n=0$, the first player plays $E^1_0=E^p$ (as the
rules of the game require, according to (e)). The second player answers
choosing an $I$--equivalence relation $E^2_0\geq E^1_0$. Now, on a side,
Player I starts to play the game of \ref{1.11} using his strategy ${\bf
St}_p$. The strategy says him to play a name $\name{\tau}_0$ of an ordinal. He
defines $q_0$ by (f) (so $q\in\qell$ is a condition stronger than $p$ and such
that $E^{q_0}=E^2_0$) and chooses a condition $p_0\geq q_0$ deciding the value
of the name $\name{\tau}_0$, say $p_0\forces\name{\tau}_0=\alpha$. He pretends
that the second player answered (in the game of \ref{1.11}) by: $p_0$,
$\bfu_0=\emptyset$, $w_0=\{\alpha\}$. Next, in the play of $GM_I(E^p)$, he
plays $E^1_1=E^{p_0}$ as declared in (e).\\
Now suppose that we are at the $(n+1)^{\rm th}$ stage of the play of
$GM_I(E^p)$, the first player has played $E^1_{n+1}$ already and on a side he
has played the play of the game \ref{1.11} as defined by (a)--(h) and ${\bf
St}_p$ (so in particular he has defined a condition $p_n$ and $E^1_{n+1}=
E^{p_n}\rest\big(\dom(E^{p_n})\setminus\bigcup\limits_{i\in\bfu_n}i/E^{p_n}
\big)$ and $\bfu_n$ is the set of the first $n$ elements of $A^{p_n}$). The
second player plays an $I$--equivalence relation $E^2_{n+1}\geq E^1_{n+1}$. 
Now the first player chooses (on a side, pretending to play in the game of
\ref{1.11}): a name $\name{\tau}_{n+1}$ given by the strategy ${\bf St}_p$, a
condition $q_{n+1}\in\qell$ determined by (f) (check that (g) is satisfied),
$\bfu_{n+1}$ as in (d) and a condition $p_{n+1}\in\qell$ satisfying (g), (h)
(the last exists by \ref{1.8}). Note that, by (g) and \ref{1.7}, the condition
$p_{n+1}$ determines a suitable set $w_{n+1}$. Thus, Player I pretends that
his opponent in the game of \ref{1.11} played $p_{n+1},\bfu_{n+1},w_{n+1}$ and
he passes to the actual game $GM_I(E^p)$. Here he plays $E^1_{n+2}$ defined
by (e). 

The strategy ${\bf St}$ described above cannot be the winning one. 
Consequently, there is a play in $GM_I(E^p)$ in which Player I uses ${\bf St}$,
but he looses. During the play he constructed a sequence $\langle
(p_n,\bfu_n,w_n):n\in\omega\rangle$ of legal moves of Player II in the game of
\ref{1.11} against the strategy ${\bf St}_p$. Let $E^q=\lim\limits_{n<\omega}
E^{p_n}$ (i.e. $\dom(E^q)=\bigcap\limits_{n<\omega}\dom(E^{p_n})$, $x\,E^q\,y$
if and only if for every large enough $n$, $x\,E^{p_n}\,y$) and let
$H^q(x^m_i)$ will be $H^{p_n}(x^m_i)$ for any large enough $n$ (it is
eventually constant). It follows from the demand (g) that $E^q$-equivalence
classes are in $I$. Moreover, $\dom(E^1_{n+1})\setminus\dom(E^2_{n+1})
\subseteq k/E^q$, where $k$ is the $(n+1)^{\rm th}$ member of $A^q$. Therefore
\[\begin{array}{l}
\omega\setminus\dom(E^q)=\omega\setminus\bigcap\limits_{n\in\omega}
\dom(E^{p_n})\subseteq\\
\omega\setminus\dom(E^{p_0})\cup \bigcup\{\dom(E^2_n)\setminus\dom(E^1_{n+1}):
n\in\omega\}\in I
  \end{array}\]
(remember, Player I lost in $GM_I(E^p)$). Now it should be clear that
$q\in\qell$ and it is stronger than every $p_n$ (even $p_n\leq^*_n q$). Hence
Player II wins the corresponding play of \ref{1.11}, showing that ${\bf St}_p$
is not a winning strategy. 
\medskip

\noindent 2)\ \ \ The same proof. \QED

\begin{proposition}
\label{1.11A} 
If in \ref{1.11} we assume $\ell = 2$ and demand $p_n \le^\otimes_n p_{n+1}$
instead $p_n \le^*_n p_{n+1}$ {\em then} Player II has a winning strategy.
\end{proposition}

\Proof Using \ref{1.8A}, the second player can find suitable conditions $p_n$
(in the game of \ref{1.11}) such that $p_n\leq^\otimes_{n+1} p_{n+1}$. But
note that the partial orders $\leq^\otimes_n$ have the fusion property, so the
sequence $\langle p_n: n<\omega\rangle$ will have an upper bound in
$\qdwa$. \QED 

\begin{remark}
\label{1.11B}
We could have used $<^\otimes_n$ also in \cite{Sh:407}.
\end{remark}

\begin{definition}[see {\cite[VI, 2.12, A-F]{Sh:f}}]
\label{1.13}
\begin{enumerate}
\item A forcing notion $\bP$ has the PP-property if:
\begin{description}
\item[($\otimes^{PP}$)]  {\em for every} $\eta\in\baire$ from $\V^{\bP}$ and a
strictly increasing $x\in \baire\cap\V$ {\em there is} a closed subtree $T
\subseteq \fseo$ such that: 
\begin{description}
\item[$(\alpha)$]  $\eta\in\lim(T)$, i.e.~$(\forall n<\omega)(\eta\restriction
n \in T)$,  
\item[$(\beta)$]  $T \cap {}^n \omega$ is finite for each $n < \omega$,
\item[$(\gamma)$] for arbitrarily large $n$ there are $k$, and $n<i(0)<j(0)<
i(1)<j(1)<\ldots<i(k)<j(k)<\omega$ and for each $\ell\le k$, there are
$m(\ell) < \omega$ and $\eta^{\ell,0},\ldots,\eta^{\ell,m(\ell)} \in T \cap
{}^{j(\ell)} \omega$ such that $j(\ell) > x(i(\ell) + m(\ell))$ and 
\[(\forall\nu\in T\cap {}^{j(k)}\omega)(\exists\ell\leq k)(\exists m\leq
m(\ell))(\eta^{\ell,m} \trianglelefteq \nu).\] 
\end{description}
\end{description}
\item We say that a forcing notion $\bP$ has the strong PP--property if
\begin{description}
\item[($\oplus^{sPP}$)]  for every function $g:\omega\longrightarrow\V$ from
$\V^{\bP}$ there exist a set $B\in\iso\cap \V$ and a sequence $\langle w_n:
n\in B\rangle\in\V$ such that for each $n\in B$
\[|w_n|\leq n\qquad\mbox{ and }\qquad g(n)\in w_n.\]
\end{description}
\end{enumerate}
\end{definition}

\begin{remark}
{\em
Of course, if a proper forcing notion has the strong PP--property then it has
the PP--property.
}
\end{remark}

\begin{conclusion}
\label{1.12}
Assume that for each $p\in\qell$ and for each countable $N\prec ({\cal H}(
\chi),{\in},{<}^*)$ such that $p,I,h\in N$, the first player has no winning
strategy in $GM^N_I(E^p)$ (e.g.~if $I$ is a maximal ideal). Then
\begin{description} 
\item[(*)]  $\qell$ is proper, $\alpha$-proper, strongly $\alpha$-proper for
every $\alpha < \omega_1$, is $\baire$-bounding and it has the PP-property,
even the strong PP--property. \QED
\end{description}
\end{conclusion}

By \cite[VI, 2.12]{Sh:f} we know

\begin{theorem}
\label{1.14}
Suppose that $\langle \bP_i,\nbQ_j:j < \alpha,i \le \alpha \rangle$ is a
countable support iteration such that 
\[\Vdash_{\bP_j} \mbox{`` }\nbQ_j\mbox{ is proper and has the PP-property}".\]
Then $\bP_\alpha$ has the PP-property. \QED
\end{theorem}

\section{NWD ultrafilters}
A subset $A$ of the set $\cQ$ of rationals is {\em nowhere dense} (NWD) if its
closure (in $\cQ$) has empty interior. Remember that the rationals are
equipped with the order topology and both ``closure'' and ``interior'' refer
to this topology. Of course, as $\cQ$ is dense in the real line, we may
consider these operations on the real line and get the same notion of nowhere
dense sets. For technical reasons, in forcing considerations we prefer to work
with $\can$ instead of the real line. So naturally we want to replace
rationals by $\fs$. But what are nowhere dense subsets of $\fs$ then? (One may
worry about the way we ``embed'' $\fs$ into $\can$.) Note that we have a
natural lexicographical ordering $<_{\lx}$ of $\fs$:

$\eta<_{\lx}\nu$\qquad if and only if

either there is $\ell<\omega$ such that $\eta\rest\ell=\nu\rest\ell$ and
$\eta(\ell)<\nu(\ell)$ 

or $\eta\conc\langle 1\rangle\trianglelefteq\nu$

or $\nu\conc\langle 0\rangle\trianglelefteq\eta$.

\noindent Clearly $(\fs,<_{\lx})$ is a linear dense order without end-points
(and consequently it is order--isomorphic to the rationals). Now, we may talk
about nowhere dense subsets of $\fs$ looking at this ordering only, but we may
relate this notion to the topology of $\can$ as well. 

\begin{proposition}
\label{equiv}
For a set $A\subseteq\fs$ the following conditions are equivalent:
\begin{enumerate}
\item $A$ is nowhere dense,
\item $(\forall \eta\in\fs)(\exists\nu\in\fs)[\eta\trianglelefteq\nu\ \&\
(\forall\rho\in\fs)(\nu\trianglelefteq\rho\ \Rightarrow\ \rho\notin A)]$,
\item the set 
\[A^*\stackrel{\rm def}{=}\{\eta\in\can: (\forall n\in\omega)(\exists\nu\in
A)(\eta\rest n\trianglelefteq\nu)\}\]
is nowhere dense (in the product topology of $\can$),
\item there is a sequence $\langle\eta_n: n<\omega\rangle$ such that for each
$n<\omega$
\begin{description}
\item[(i)$_n$]   $\eta_n:[n,\ell_n)\longrightarrow 2$ for some $\ell_n>n$ and
\item[(ii)$_n$]  $(\forall\rho\in A)(\eta_n\not\subseteq\rho)$,
\end{description}
\item there is a sequence $\langle\eta_n: n<\omega\rangle$ such that for each
$n<\omega$ condition {\bf (i)}$_n$ (see above) holds and
\begin{description}
\item[(ii)$^*_n$]  $(\forall\nu\in {}^{\textstyle n}2)(\{\rho\in\fs: \nu\cup
\eta_n\trianglelefteq\rho\}\cap A=\emptyset)$,
\end{description}
\item there are $B\in\iso$ and $\langle\eta_n: n\in B\rangle$ such that for
each $n\in B$ the conditions {\bf (i)}$_n$, {\bf (ii)}$_n$ above are
satisfied. 
\end{enumerate}
\end{proposition}

\Proof $1.\ \Rightarrow\ 2.$\qquad Suppose $A\subseteq\fs$ is nowhere dense
but for some sequence $\eta\in\fs$, for every $\nu\in\fs$ extending $\eta$
there is $\rho\in A$ such that $\nu\trianglelefteq\rho$. Look at the interval
$(\eta\conc\langle 0\rangle,\eta\conc\langle 1\rangle)_{<_{\lx}}$ (of
$(\fs,<_{\lx})$). We claim that $A$ is dense in this interval. Why? Suppose
\[\eta\conc\langle 0\rangle\leq_{\lx}\eta^*_0<_{\lx}\eta^*_1 \leq_{\lx}
\eta\conc\langle 1\rangle.\]
Assume $\lh(\eta^*_0)\leq\lh(\eta^*_1)$. Take $\nu\stackrel{\rm def}{=}
\eta^*_1\conc\langle 0\rangle$. By the definition of the order $<_{\lx}$ we
have then
\[\eta^*_0<_{\lx}\nu\conc\langle 0\rangle<_{\lx}\nu\conc\langle 1\rangle
<_{\lx}\eta^*_1\qquad\mbox{ and }\qquad\eta\vartriangleleft\nu.\]
By our assumption we find $\rho\in A$ such that $\nu\conc\langle0,1\rangle
\trianglelefteq\rho$. Then
\[\nu\conc\langle 0\rangle<_{\lx}\rho<_{\lx}\nu\conc\langle 1\rangle
\qquad\mbox{ and hence }\quad \rho\in (\eta^*_0,\eta^*_1)_{<_{\lx}}.\]
Similarly if $\lh(\eta^*_1)\leq\lh(\eta^*_0)$.

$2.\ \Rightarrow\ 3.$\qquad Should be clear if you remember that sets
\[[\nu]\stackrel{\rm def}{=} \{\eta\in\can: \nu\vartriangleleft\eta\}\qquad
\qquad\mbox{ (for $\nu\in\fs$)}\]
constitute the basis of the topology of $\can$.

$3.\ \Rightarrow\ 4.$\qquad Suppose $A^*$ is nowhere dense in $\can$. Let
$n<\omega$. Considering all elements of $2^{\textstyle n}$ build
(e.g. inductively) a function $\eta^*_n: [n,\ell^*_n)\longrightarrow 2$ such
that $n<\ell^*_n$ and 
\[(\forall\nu\in 2^{\textstyle n})([\nu\conc\eta^*_n]\cap A^*=\emptyset).\]
This means that for each $\nu\in 2^{\textstyle n}$ the set $\{\rho\in A:
\nu\conc\eta^*_n\trianglelefteq\rho\}$ is finite (otherwise use K\"onig lemma
to construct an element of $A^*$ in $[\nu\conc\eta^*_n]$). Taking sufficiently
large $\ell_n>\ell^*_n$ and extending $\eta^*_n$ to $\eta_n$ with domain
$[n,\ell_n)$ we get that $(\forall\rho\in A)(\eta_n\not\subseteq\rho)$ (as
required). 

$4.\ \Rightarrow\ 5.\ \Rightarrow\ 6.$\qquad Read the conditions.

$6.\ \Rightarrow\ 1.$\qquad Let $B$, $\langle\eta_n:n\in B\rangle$ be as in
6. Suppose $\nu_0,\nu_1\in\fs$, $\nu_0<_{\lx}\nu_1$. Assume
$\lh(\nu_0)\leq\lh(\nu_1)=m$. Take any $n\in B\setminus (m+1)$ and let
$\nu=\nu_1\conc\langle\underbrace{0,\ldots,0}_{n-m}\rangle\conc\eta_n$. We
know that no element of $A$ extends $\nu$. But this implies that the interval
$(\nu\conc\langle 0\rangle, \nu\conc\langle 1\rangle)_{<_{\lx}}$ is disjoint
from $A$ (and is contained in the interval
$(\nu_0,\nu_1)_{<_{\lx}}$). Similarly if $\lh(\nu_1)\leq\lh(\nu_0)$. \QED

\begin{lemma}
\label{2.3}
Let $n,k^*<\omega$. Assume that $\bar{\nu}^k=\langle\nu^k_i:n\le i<i_k
\rangle$ for $k<k^*<\omega$, $n\le i_k<\omega$, $\nu^k_i\in\bigcup\limits_{j
\ge i}{}^{[i,j)}2$ and $w_k \subseteq [n,i_k)$, $|w_k| \ge k^*$ and: 
\[\mbox{if }k< k^*,\ m_1< m_2\mbox{ are in $w_k$ then }\max\dom(\nu^k_{m_1})<
m_2.\]
Lastly let 
\[i(*)=\max\{\sup\dom(\nu^k_i)+1:k< k^* \mbox{ and } i \in (n,i_k)\}.\] 
{\em Then} we can find $\rho \in {}^{[n,i(*))} 2$ such that: 
\[(\forall k < k^*)(\exists i \in w_k)(\nu^k_i \subseteq \rho).\]
\end{lemma}

\Proof  By induction on $k^*$ (for all possible other parameters). For $k^*
=0,1$ it is  trivial.\\
Let $n^0_k=\min(w_k)$ and $n^1_k=\min(w_k\setminus(n^0_k+1))$. Let $\ell<k^*$
be with minimal $n^1_\ell$. Apply the induction hypothesis with
$n^1_\ell$, $\bar{\nu}^k=\langle\bar{\nu}^k_i:n^1_\ell\le i<i_k\rangle$ for $k
<k^*$, $k \ne \ell$ and $\langle w_k\setminus n^1_\ell:k<k^*,k\neq\ell\rangle$
here standing for $n$, $\bar{\nu}^k$ for $k< k^*$, $\langle w_k:k<k^*\rangle$
there and get $\rho_1\in {}^{[n^1_\ell,i(*))}2$. Note that $w_k\setminus
n^1_\ell\supseteq w_k \setminus n^1_k$ has at least $|w_k|-1$ elements. Let
$\rho \in {}^{[n,i(*))}2$ be such that $\rho_1 \subseteq \rho$ and
$\bar{\nu}^\ell_{n^0_\ell} \subseteq \rho$. \QED

\begin{proposition}
\label{2.3A}
Assume that $\bR$ is a proper forcing notion with the PP-property. Then
\begin{description}
\item[($\oplus^{\nwd}$)] for every nowhere dense set $A\subseteq\fs$ in
$\V^{\bR}$ there is a nowhere dense set $A^*\subseteq\fs$ in $\V$ such that
$A\subseteq A^*$.
\end{description}
\end{proposition}

\Proof  Let $A\in\V^{\bR}$ be a nowhere dense subset of $\fs$. Thus, in
$\V^{\bR}$, we can, for each $n<\omega$, choose $\nu_n\in\bigcup\limits_{\ell
\ge n} {}^{[n,\ell)}2$ such that:
\[(\forall\nu\in {}^n 2)(\forall \rho\in\fs)(\nu\conc \nu_n\trianglelefteq
\rho\ \Rightarrow\ \rho \notin A).\]
So $\langle \nu_n:n<\omega\rangle\in \V^{\bR}$ is well defined. Next for each
$n$ we choose an integer $\ell_n \in (n,\omega)$, a sequence $\eta_n\in
{}^{[n,\ell_n)}2$ and a set $w_n \subseteq [n,\ell_n)$ such that:
\begin{itemize}
\item $|w_n| > n$,
\item $(\forall m\in w_n)(\nu_m\subseteq\eta_n)$, so in particular $(\forall m
\in w_n)(\max\dom(\nu_m)<\ell_n)$, and 
\item for any $m_1<m_2$ from $w_n$ we have $\max\dom(\nu_{m_1})<m_2$.
\end{itemize}
So $\bar{w}=\langle w_n:n<\omega\rangle,\,\bar{\eta}=\langle\eta_n:n<\omega
\rangle\in \V^{\bR}$ are well defined. 

Since $\bR$ has the PP-property it is $\baire$-bounding, and hence there is a
strictly increasing $x\in\baire\cap\V$ such that $(\forall n\in\omega)(\ell_n
<x(n))$.  Applying the PP-property of $\bR$ to $x$ and the function $n\mapsto
(\eta_n,w_n)$ we can find $\langle\langle V^n_\ell:\ell \le k_n \rangle:n<
\omega\rangle$ in $\V$ and $\langle\langle(i_\ell(n),j_\ell(n)):\ell\le k_n
\rangle:n<\omega\rangle$ in $\V$
such that: 
\begin{description}
\item[(a)]   $i_0(n)<j_0(n)<i_1(n)<j_1(n)<\ldots<i_{k_n}(n)<j_{k_n}(n)$,
\item[(b)]   $j_{k_n}(n)<i_0(n+1)$ for $n<\omega$,
\item[(c)]   $x(i_\ell(n))<j_\ell(n)$,
\item[(d)]   $V^n_\ell\subseteq\{(\eta,w):\eta \in {}^{[i_\ell(n),j_\ell(n))}
2$ and $w \subseteq [i_\ell(n),j_\ell(n)),\ |w|>i_\ell(n)\}$ for $\ell\le
k_n$, $n<\omega$, 
\item[(e)]   $|V^n_\ell|\le i_\ell(n)$,
\item[(f)]   for every $n<\omega$, for some $\ell\le k_n$ and $(\eta,w)\in
V^n_\ell$ we have $w = w_{i_\ell(n)}$, $\eta_{i_\ell(n)} \subseteq \eta$. 
\end{description}
[Note that $i_\ell(n)$ corresponds to $i(\ell)+m(\ell)$ in definition
\ref{1.13}(1), so we do not have $m_\ell(n)$ here.] Working in $\V$, by
\ref{2.3}, for each $n<\omega$, $\ell\le k_n$ there is 
$\rho^n_\ell\in {}^{[i_\ell(n),j_\ell(n))}2$ such that: 
\[(\forall(\eta,w)\in V^n_\ell)(\exists m_1,m_2\in w)(m_2=\min(w\setminus (m_1
+ 1))\ \&\ \eta\restriction [m_1,m_2)\subseteq \rho^n_\ell).\]
Let $\rho_n\in {}^{[i_0(n),i_0(n+1))}2$ be such that $\ell\le k_n\ \
\Rightarrow\ \ \rho^n_\ell\subseteq \rho_n$. As we have worked in $\V$,
$\langle \rho_n:n<\omega\rangle \in \V$.  Let 
\[A^* = \{\rho\in\fs: \neg(\exists n\in\omega)(\rho_n \subseteq \rho)\}.\]
Clearly $A^* \in\V$ is as required. \QED
\medskip

Let us recall definition \ref{0.1} reformulating it slightly for technical
purposes. (Of course, the two definitions are equivalent; see the discussion
at the beginning of this section.)

\begin{definition}
\label{2.2}
We say that a non-principal ultrafilter $\D$ on $\omega$ is an NWD-ultrafilter
{\em if} for any sequence $\langle\eta_n:n<\omega\rangle\subseteq\fs$ for some
$A\in\D$ the set $\{\eta_n:n\in A\}$ is nowhere dense in $\fs$.
\end{definition}

\begin{lemma}
\label{2.1}
Let $\D$ be a non-principal ultrafilter on $\omega$ and $I$ be the dual ideal
(and $h:\omega\longrightarrow\omega$ non-decreasing $\lim\limits_{n\to\infty}
h(n)=\infty$).  Then:
\begin{enumerate}
\item  in $\V^{\qjeden}$ we cannot extend $\D$ to an NWD-ultrafilter.
\item  If $\nbQ$ is a $\qjeden$-name of a proper forcing notion with the
PP--property, {\em then} also in $\V^{\qjeden * \nbQ}$ we cannot extend $\D$
to an NWD-ultrafilter. 
\end{enumerate}
\end{lemma}

\Proof  1)\ \ \ Let $\name{\bar{\eta}}=\langle\name{\eta}_n: n<\omega\rangle$
be the name defined in \ref{1.4}, but now we interpret the value $-1$ as
$0$. So $\Vdash$``$\name{\eta}_n\in {}^{h(n)} 2$'' (for each
$n<\omega$). Clearly it is enough to show that
\[\begin{array}{ll}
(*)\qquad\forces_{\qjeden}&\mbox{`` if $X\subseteq\omega$ and the set
$\{\name{\eta}_n: n\in X\}$ is nowhere dense}\\
&\mbox{ \ then there is $Y\in\D$ disjoint from $X$''.}
  \end{array}\]
So suppose that $\name{\tau}$ is a $\qjeden$-name for a subset of $\omega$ and
a condition $p^*\in\qjeden$ forces that $\{\name{\eta}_n: n\in\name{\tau}\}$
is nowhere dense. By \ref{equiv}, for some $\qjeden$-names $\name{\bar{\nu}}=
\langle\name{\nu}_m:m<\omega\rangle$ we have 
\[p^* \Vdash\mbox{``}\name{\nu}_m\in\bigcup_{\ell\ge m}{}^{[m,\ell)}2\mbox{
and for every $m<\omega$ for no $n\in\name{\tau}$ we have $\name{\nu}_m
\subseteq \name{\eta}_n$''.}\]
By \ref{1.11} (or actually by its proof) without loss of generality: 
\begin{quotation}
\noindent for every $n\in A^{p^*}$, for some $k_n\in (n,\min(A^{p^*}\setminus
(n+1)))$, for every $f:\{x^m_j:m\in A^{p^*}\cap (n+1)\mbox{ and } j<h(m)\}
\longrightarrow \{-1,1\}$, the condition $p^{*^{[f]}}$ forces a value to
$\name{\tau}\cap k_n$, and $\name{\tau}\cap k_n\setminus n\ne\emptyset$.   
\end{quotation}
[Why? Give a strategy to Player I in the game there for $p^*$ trying to force
the needed information, so for some such play Player II wins and replaces
$p^*$ by $q$ from there.]\\
Again by \ref{1.11} we may assume that
\begin{quotation}
\noindent for every $f:\{x^m_j:j<h(m)\mbox{ and } m\in A^{p^*}\cap (n+1)\}
\longrightarrow\{-1,1\}$, $n \in A^{p^*}$, for some $\bar{\nu}^f$ we have 
\[p^{*^{[f]}}\Vdash\mbox{``}\bar{\nu}^f\mbox{ is an initial segment of }
\name{\bar{\nu}}\mbox{ and }\lh(\bar{\nu}^f)= n+1\mbox{ ''.}\]
\end{quotation}
For $n\in A^{p^*}$ and $f:\{x^m_j:j<h(m)\mbox{ and }m\in A^{p^*}\cap
(n+1)\}\longrightarrow \{-1,1\}$ and $k\in A^{p^*}\setminus (n+1)$ let:
\begin{description}
\item[(a)]  $f^{[k,p^*]}$ be the function with domain $\{x^m_j:j< h(m)\mbox{
and } m \in A^{p^*} \cap (k+1)\}$ extending $f$ that is constantly 1 on
$\dom(f^{[k,p^*]})\setminus\dom(f)$,
\item[(b)]  $\bar{\rho}^f$ be an $\omega$-sequence $\langle\rho^f_\ell:\ell<
\omega\rangle$ such that for each $k\in A^{p^*}\setminus(n+1)$ we have
$\bar{\rho}^f\restriction (k+1)=\bar{\nu}^{f^{[k,p^*]}}\restriction (k+1)$. 
\end{description}
Now, for every $n \in A^{p^*}$, we can find $\rho^*_n\in\fs$ such that for
every function
\[f:\{x^m_j:j < h(m) \mbox{ and } m \in A^{p^*} \cap (n+1)\} \longrightarrow
\{-1,1\}\]
for some $\ell(f) \in (h(n),\omega)$ we have $\rho^f_{\ell(f)}\subseteq
\rho^*_n$ (so $\ell(f)<\ell g(\rho^*_n)$).\\
\relax [Why? Let $\{f_j:j < j(*)\}$ list the possible $f$'s, and we chose by
induction on $j \le j(*)$, $\rho^j\in\fs$ such that $\rho^j\vartriangleleft
\rho^{j+1}$, and $\rho^{j+1}$ satisfies the requirement on $f_j$, e.g. $\rho_0
= \langle\underbrace{0,\ldots,0}_{h(n)}\rangle$, $\rho^{j+1} = \rho^j\conc
\rho^{f_j}_{\ell g(\rho^j)}$]. 

Now choose by induction on $\zeta<\omega$, $n_\zeta\in A^{p^*}$ such that
$n_\zeta<n_{\zeta+1}$, and $\lh(\rho^*_{n_\zeta})<h(n_{\zeta +1})$. Without
loss of generality $\bigcup\limits_{\zeta<\omega}(n_\zeta/E^{p^*})\in I$. Then
\begin{quotation}
\noindent either\ \ $\bigcup\{n/E^{p^*}:n\in A^{p^*} \mbox{ and }(\exists\zeta
<\omega)(n_{2 \zeta}<n<n_{2\zeta+1})\}\in\D$

\noindent or\qquad $\bigcup\{n/E^{p^*}:n\in A^{p^*} \mbox{ and }(\exists\zeta<
\omega)(n_{2 \zeta+1}<n<n_{2\zeta+2})\}\in\D$,
\end{quotation}
so by renaming the latter holds. (Again, it suffices that the ideal $I$ is
such that the quotient algebra ${\cal P}(\omega)/I$ satisfies the c.c.c.)
Lastly we define a condition $r\in\qjeden$: 
\[\dom(E^r)=\bigcup_{\zeta<\omega}n_{2\zeta}/E^{p^*}\cup \bigcup\{n/E^{p^*}\!:
n\!\in\! A^{p^*}\mbox{ and }(\exists\zeta\!<\!\omega)(n_{2\zeta+1}<n<n_{2
\zeta+2})\},\] 
\[n_{2 \zeta}/E^r=(n_{2 \zeta}/E^{p^*})\cup \bigcup\{m/E^{p^*}:m\in A^{p^*}
\cap (n_{2 \zeta+1},n_{2 \zeta+ 2})\}\]
(note that this defines correctly an $I$--equivalence relation $E^r$), $A^r=
\{n_{2 \zeta}:\zeta<\omega\}$. The function $H^r$ is defined by cases
(interpreting the value $0$ as $-1$, where appears):
\[\begin{array}{lcr}
H^r(x^m_j)= H^{p^*}(x^m_j) & \mbox{if}
&m\in(\omega\setminus\dom(E^{p^*}))\mbox{ and } j<h(m),\\
H^r(x^m_j)=H^{p^*}(x^m_j) & \mbox{if} & m\in\dom(E^{p^*})\mbox{ and } j\in
[h(\min(m/E^{p^*})),h(m))\\
H^r(x^m_j)=1 & \mbox{if} & m\in\dom(E^{p^*})\mbox{ and }\min(m/E^{p^*})\in
(n_{2 \zeta},n_{2 \zeta +1}]\\
&  &\mbox{and } j<h(\min(m/E^{p^*}))\\
H^r(x^m_j)=\rho^*_{n_{2 \zeta}}(j) & \mbox{if} & m\in\dom(E^{p^*})\mbox{ and }
\min(m/E^{p^*})\in (n_{2 \zeta +1},n_{2 \zeta + 2})\\
&  &\mbox{and } j\in\dom(\rho^*_{n_{2 \zeta}})\mbox{ and }j\ge h(n_{2\zeta})\\
H^r(x^m_j)=1 & \mbox{ } & \mbox{otherwise (but }x^m_j \in \dom(H^r)).
  \end{array}\]
Now check that $p^* \le r \in\qjeden$ and for each $n\in\dom(E^r)\setminus
\bigcup\limits_{\zeta<\omega} n_{2\zeta}/E^{p^*}$:
\[r\Vdash\mbox{`` }\name{\eta}_n\mbox{ violates the property of }\name{\bar{
\nu}}\mbox{ and hence }n\notin\name{\tau}\mbox{''.}\]
As $\dom(E^r)\setminus\bigcup\limits_{\zeta<\omega} n_{2\zeta}/E^{p^*}\in\D$
we have finished.  
\medskip

\noindent 2)\ \ \ Should be clear by (*) of the proof of \ref{2.1}(1) and
\ref{2.3A}.\\ 
However we will give an alternative proof of \ref{2.1}(2). We start as in
the proof of \ref{2.1}(1): suppose some $(p^*,\name{r}^*)\in\qjeden * \nbQ$
forces ``$\name{F}$ is an NWD-ultrafilter on $\omega$ extending $\D$''.  As
$\Vdash$``$\name{\eta}_n[\name{G}_{\qjeden}]\in {}^{h(n)}2$'', for some
$(\qjeden * \nbQ)$-name $\name{\tau}$ for a subset of $\omega$ 
\[(p^*,\name{r}^*)\Vdash\mbox{`` }\name{\tau}\in\name{F}\mbox{ and }(\forall
\eta \in\fs)(\exists \nu\in\fs)(\eta\trianglelefteq\nu\ \&\ (\forall n\in 
\name{\tau})(\neg\nu\trianglelefteq\name{\eta}_n))\mbox{ ''.}\]
So for some $\qjeden*\nbQ$-name $\name{\bar{\nu}}=\langle\name{\nu}_n:n< 
\omega\rangle$
\[(p^*,\name{r}^*)\Vdash\mbox{`` }\name{\nu}_\ell\in\bigcup_{j\in [\ell,
\omega)} {}^{[\ell,j)}2\mbox{ and for no }n\in \name{\tau}\mbox{ we have
}\name{\nu}_\ell\subseteq \name{\eta}_n\mbox{''.}\]
So for some $\qjeden*\nbQ$--names $\name{d}_\ell$, $\name{w}_\ell$
\[\begin{array}{ll}
(p^*,\name{r}^*)\Vdash&\mbox{`` }\omega>\name{d}_\ell>\ell,\ \name{w}_\ell
\subseteq [\ell,\name{d}_\ell),\ |\name{w}_\ell|>(4\cdot\prod_{s\le
n}h(s))!\mbox{ and}\\
&[m_1 < m_2 \mbox{ in } \name{w}_\ell\quad\Rightarrow\quad \max\dom(\name{
\nu}_{m_1})<m_2]\mbox{''.}
  \end{array}\]
Let $p^* \in G_{\qjeden}\subseteq \qjeden$ and $G_{\qjeden}$ generic over
$\V$. Now in $\V[G_{\qjeden}]$, the forcing notion $\nbQ[G_{\qjeden}]$ is
$\baire$-bounding (this follows from the PP-property) and also $\qjeden$ is 
$\baire$-bounding. Hence for some $r'\in \nbQ[G_{\qjeden}]$ and strictly
increasing $x \in \baire\cap\V$ we have: 
\[r'\Vdash_{\nbQ[G_{\qjeden}]}\mbox{`` }\name{d}_n<x(n)\mbox{ and } m\in \name{w}_n\ \Rightarrow\ \dom(\name{\nu}_m)\subseteq [0,x(n))\mbox{''.}\]
In $\V[G_{\qjeden}]$, by the property of $\nbQ$, there are $r^{**}$, $r'\le
r^{**} \in \nbQ[G_{\qjeden}]$ and a sequence $\langle\langle
i_\ell(n),j_\ell(n)):\ell \le k_n \rangle:n < \omega\rangle$ such that 
\[i_0(n) < j_0(n) < i_1(n) < j_1(n) < \ldots < j_{k_n}(n) < i_\ell
(n+1),j_\ell(n) > x(i_\ell(n))\]
and there are $\bar{\nu}^*_{n,\ell,t} = \langle \nu^*_{n,\ell,t,j}:
j \in [i_\ell(n),j_\ell(n)) \rangle$ for $t < i_\ell(n),\ell \le k_n$ and
$\bar{w}^*_{n,\ell,t} = \langle w^*_{n,\ell,t,j}:j \in [i_\ell(n),i_{\ell + 1}
(n))$ for $t < i_\ell(n),\ell \le k_n \rangle$ such that 
\[\begin{array}{ll}
r^{**} \Vdash_{\nbQ}&\mbox{``}\langle\name{\nu}_{i_\ell(n)+j}:j\in [i_\ell(n),
j_\ell(n))\rangle \mbox{ is }\bar{\nu}^*_{n,\ell,t}\mbox{ and}\\
&\langle \name{w}_{i_\ell(n) + j}:j \in [i_\ell(n):j_\ell(n))\rangle\mbox{ is
}\bar{w}^*_{n,\ell,t}\mbox{ for some }t < i_\ell(n)\mbox{''.}
  \end{array}\]
Back in $\V$ we have a $\qjeden$-name $\name{r}^{**}$ and $\langle\langle(
\name{i}_\ell(n),\name{j}_\ell(n)):\ell\le \name{k}_n\rangle:n<\omega\rangle$
and $\langle\langle \name{\bar{\nu}}^*_{n,\ell,t}:t<i_\ell(n)\rangle:\ell<
\name{k}_n,n<\omega\rangle$ and $\langle\langle\bar{w}^*_{n,\ell,t}:t<
i_\ell(n)\rangle: \ell<\name{k}_n,n<\omega\rangle$ are forced (by $p^*$) to be
as above. 

By \ref{1.11}, increasing $p^*$, we get
\begin{quotation}
\noindent for every $f:\{x^n_i:i < h(m),m \in A^{p^*} \cap (n+1)\}
\longrightarrow\{-1,1\}$, $n \in A^{p^*}$, the condition $p^{*^{[f]}}$ forces
a value to  
\[\begin{array}{l}
\langle\langle(\name{i}_\ell(m),\name{j}_\ell(m)):\ell \le\name{k}_m\rangle:
m\le n\rangle,\\
\langle\name{\bar{\nu}}^*_{n,\ell,t}:t<\name{i}_\ell(n),\ell\le\name{k}_n
\rangle,\\
\langle\name{\bar{w}}^*_{m,\ell,t}:t <\name{i}_\ell(n),\ell<\name{k}_n\rangle
  \end{array}\]
moreover, without loss of generality 
\[n \in A^{p^*}\quad \Rightarrow j_{\name{k}_n}(n)<\min(A^{p^*}\setminus
(n+1)).\]
\end{quotation}
Now by \ref{2.3}, without loss of generality for each $n \in A^{p^*}$ we can
find a function $\rho_n$ from $[n,\min(A^{p^*}\setminus (n+1))]$ to $\{-1,1\}$
such that: 
\begin{quotation}
\noindent if $f:\{x^m_i:i < h(m),m \in A^{p^*}\cap(n+1)\}\longrightarrow\{-1,
1\}$, $n \in A^{p^*}$

\noindent then $(p^{*^{[f]}},\name{r}^{**})$ forces that $\rho_n$ extends some
$\name{\nu}_\ell$.
\end{quotation}
Now we continue as in the proof of \ref{2.1}(1). \QED

\section{The consistency proof}

\begin{theorem}
\label{2.4}
Assume CH and $\diamondsuit_{\{\gamma<\omega_2: \cf(\gamma)=\omega_1\}}$.

\noindent Then there is an $\aleph_2$--cc proper forcing notion $\bP$ of
cardinality $\aleph_2$ such that
\[\forces_{\bP}\mbox{`` there are no NWD--ultrafilters on $\omega$ ''.}\]
\end{theorem}

\Proof Define a countable support iteration $\langle \bP_i,\nbQ_j:i\leq
\omega_2, j<\omega_2\rangle$ of proper forcing notions and sequences
$\langle\nD_i: i<\omega_2\rangle$ and $\langle\name{\bar{\eta}}^i: i<\omega_2
\rangle$ such that for each $i<\omega_2$: 
\begin{enumerate}
\item $\nD_i$ is a $\bP_i$--name for a non--principal ultrafilter on $\omega$, 
\item $\nbQ_i$ is a $\bP_i$--name for a proper forcing notion of size
$\aleph_1$ with the PP--property,
\item $\name{\bar{\eta}}^i$ is a $\bP_i*\nbQ_i$--name for a function from
$\omega$ to $\fs$, 
\item $\forces_{\bP_i*\nbQ_i}$ ``if $X\subseteq\omega$ and the set
$\{\name{\eta}^i_n: n\in X\}\subseteq\fs$ is nowhere dense then there is
$Y\in\nD_i$ disjoint from $X$'',  
\item if $\nD$ is a $\bP_{\omega_2}$--name for an ultrafilter on $\omega$ then
the set 
\[\{i<\omega_2:\ \cf(i)=\omega_1\quad\&\quad\nD_i =\nD\rest {\cal
P}(\omega)^{\V^{\bP_i}}\}\] 
is stationary.
\end{enumerate}
Let us first argue that if we succeed with the construction then, in
$\V^{\bP_{\omega_2}}$, we will have
\[2^{\aleph_0}=\aleph_2\ \ + \ \mbox{``there is no NWD-ultrafilter on 
}\omega\mbox{''.}\] 
Why? As each $\nbQ_i$ is (a name) for a proper forcing notion of size
$\aleph_1$, the limit $\bP_{\omega_2}$ is a proper forcing notion with a dense
subset of size $\aleph_2$ and satisfying the $\aleph_2$--cc. Since
$\bP_{\omega_2}$ is proper, each subset of $\omega$ (in $\V^{\bP_{\omega_2}}$)
has a canonical countable name (i.e. a name which is a sequence of countable
antichains; every condition in the $n^{\rm th}$ antichain decides if the
integer $n$ is in the set or not; of course we do not require that the
antichains are maximal). Hence $\forces_{\bP_{\omega_2}}2^{\aleph_0}\leq
\aleph_2$ (remember that we have assumed $\V\models$CH). Moreover, by
\ref{1.14} + \ref{2.3A} we know that $\bP_{\omega_2}$ satisfies
$(\oplus^{\nwd})$ of \ref{2.3A}, i.e. 
\[\begin{array}{ll}
\forces_{\bP_{\omega_2}}&\mbox{``each nowhere dense subset of $\fs$ can be
covered}\\
\ &\mbox{\ \ by a nowhere dense subset of $\fs$ from $\V$''}.
  \end{array}
\]
Now suppose that $\nD$ is a $\bP_{\omega_2}$--name for an ultrafilter on
$\omega$. By the fifth requirement, we find $i<\omega_2$ such that $\nD_i=
\nD\rest{\cal P}(\omega)^{\V^{\bP_i}}$ (and $\cf(i)=\omega_1$). Since
$\bP_{\omega_2}$ satisfies $(\oplus^{\nwd})$, we have
\begin{description}
\item[$\forces_{\bP_{\omega_2}}$] ``if $X\subseteq\omega$ and the set
$\{\name{\eta}^i_n: n\in X\}\subseteq\fs$ is nowhere dense then there is an
element of $\nD\rest {\cal P}(\omega)^{\V^{\bP_i}}$ disjoint from $X$''
\end{description}
[Why? Cover $\{\name{\eta}^i_n: n\in X\}$ by a nowhere dense set $A\subseteq
\fs$ from $\V$ and look at the set $X_0=\{n\in\omega: \name{\eta}^i_n\in
A\}$. Clearly $X_0\in\V^{\bP_i*\nbQ_i}$ and $X\subseteq X_0$. Applying the
fourth clause to $X_0$ we find $Y\in\nD_i=\nD\rest{\cal
P}(\omega)^{\V^{\bP_i}}$ such that $Y\cap X_0=\emptyset$. Then $Y\cap
X=\emptyset$ too.]\\
But this means that, in $\V^{\bP_{\omega_2}}$, the function
$\name{\bar{\eta}}^i$ exemplifies that $\nD$ is not an NWD ultrafilter
(remember $\nD\rest{\cal P}(\omega)^{\V^{\bP_i}}\subseteq\nD$). Moreover, as
CH implies the existence of NWD-ultrafilters, we conclude that actually
$\forces_{\bP_{\omega_2}} 2^{\aleph_0}=\aleph_2$. 

Let us describe how one can carry out the construction. Each $\nbQ_i$ will be
$\bQ^1_{\name{I}_i,h}$ for some increasing function $h\in\baire$
(e.g. $h(n)=n$) and a ($\bP_i$--name for a) maximal non--principal ideal
$\name{I}_i$ on $\omega$. By \ref{2.2}, \ref{1.12} we know that $\nbQ^1_{
\name{I}_i,h}$ satisfies the demands (2)--(4) for the ultrafilter $\nD_i$ dual
to $\name{I}_i$ and the function $\name{\bar{\eta}}^i$ as in the proof of
\ref{2.2}. Thus, what we have to do is to say what are the names $\nD_i$. To
choose them we will use the assumption of $\diamondsuit_{\{\gamma<\omega_2:
\cf(\gamma)=\omega_1\}}$. In the process of building the iteration we choose
an enumeration $\langle(p_i,\name{\tau}_i):i<\omega_2\rangle$ of all pairs 
$(p,\name{\tau})$ such that $p$ is a condition in $\bP_{\omega_2}$ (in its
standard dense subset of size $\aleph_2$) and $\name{\tau}$ is a canonical
(countable) $\bP_{\omega_2}$--name for a subset of $\omega$. We require that
$p_i\in\bP_i$ and $\name{\tau}_i$ is a $\bP_i$--name (of course, it is done by
a classical bookkeeping argument). Note that each subset of ${\cal P}(\omega)$
from $\V^{\bP_{\omega_2}}$ has a name which may be interpreted as a subset $X$
of $\omega_2$: if $i\in X$ then $p_i$ forces that $\name{\tau}_i$ is in our
set. Now we may describe how we choose the names $\nD_i$. By
$\diamondsuit_{\{\gamma<\omega_2:\cf(\gamma)= \omega_1\}}$ we have a sequence
$\langle X_i: i<\omega_2\ \&\ \cf(i)=\omega_1\rangle$ such that
\begin{description}
\item[(i)]  $X_i\subseteq i$ for each $i\in\omega_2$, $\cf(i)=\omega_1$,
\item[(ii)] if $X\subseteq\omega_2$ then the set
\[\{i\in\omega_2: \cf(i)=\omega_1\ \&\ X_i= X\cap i\}\]
is stationary.
\end{description}
Arriving at stage $i<\omega_2$, $\cf(i)=\omega_1$ we look at the set
$X_i$. We ask if it codes a $\bP_i$--name for an ultrafilter on $\omega$
(i.e. we look at $\{(p_\alpha,\name{\tau}_\alpha): \alpha\in X_i\}$ which may
be interpreted as a $\bP_i$--name for a subset of ${\cal P}(\omega)$). If yes,
then we take this name as $\nD_i$. In all remaining cases we take whatever we
wish, we may even not define the name $\name{\bar{\eta}}^i$ (note: this leaves
us a lot of freedom and one may use this to get some additional properties of
the final model). So why we may be sure that the fifth requirement is
satisfied? Suppose that we have a $\bP_{\omega_2}$--name for an ultrafilter on
$\omega$. This name can be thought of as a subset $X$ of $\omega_2$. If
$i<\omega_2$ is sufficiently closed then $X\cap i$ is a $\bP_i$--name for an
ultrafilter on $\omega$ which is the restriction of $\nD$ to $\V^{\bP_i}$. So
we have a club $C\subseteq\omega_2$ such that for each $i\in C$, if
$\cf(i)=\omega_1$ the $X\cap i$ is of this type. By {\bf (ii)} the set  
\[S\stackrel{\rm def}{=}\{i<\omega_2: i\in C\ \&\ \cf(i)=\omega_1\ \&\
X_i=X\cap i\}\]
is stationary. But easily, for each $i\in S$, the name $\nD_i$ has been chosen
in such a way that $\nD_i=\nD\rest{\cal P}(\omega)$, so we are done. \QED
\bigskip

\noindent We note that this implies that there is also no ultrafilter with
property $M$.  This was asked by Benedikt in \cite{Bn}. 

\begin{definition}
\label{2.5}
A non-principal ultrafilter $\D$ on $\omega$ has the $M$-property (or
property $M$) if:
\begin{quotation}
\noindent if for some real $\varepsilon > 0$, for $n < \omega$ we have
a tree $T_n\subseteq\fs$ such that $\mu(\lim(T_n))\ge\varepsilon$ 

\noindent then $(\exists A\in\D)(\bigcap\limits_{n\in A}\lim(T_n)\ne
\emptyset)$
\end{quotation}
(where $\mu$ stands for the Lebesgue measure on $\can$).
\end{definition}

\begin{proposition}
\label{2.6}
If a non-principal ultrafilter $\D$ on $\omega$ is not NWD, then $\D$ does not
have the property $M$.
\end{proposition}

\Proof Let
\[S^\varepsilon_\ell = \big\{T\cap {}^{\ell \ge}2:T \subseteq\fs,\ T\mbox{ a
tree not containing a cone, } \mu(\lim(T)) > \varepsilon \big\}\]
(note that $S^\varepsilon_\ell$ is a set of trees not a set of nodes) and let 
$S^\varepsilon =\bigcup\limits_\ell S^\varepsilon_\ell$.

\noindent Now let $t_1\prec t_2$ if:\quad $t_1\in S^\varepsilon_{\ell_1}$,
$t_2\in S^\varepsilon_{\ell_2}$, $\ell_1<\ell_2$ and $t_1=t_2\cap {}^{\ell_1
\ge}2$.  So $S^\varepsilon$ is a tree with $\omega$ levels, each level is
finite. As $\D$ is not NWD, we can find $\eta^\varepsilon_n\in\lim(
S^\varepsilon)$ for $n<\omega$ such that: 
\begin{quotation}
if $A \in \D$  then $\{ \eta^\varepsilon_n:n \in A\}$ is somewhere dense.
\end{quotation}
Now let $T^\varepsilon_n \subseteq \fs$ be a tree such that $\langle
T^\varepsilon_n \cap {}^{\ell \ge}2:\ell<\omega\rangle=\eta^\varepsilon_n$.
We claim that: 
\begin{quotation}
$\langle T^\varepsilon_n:n<\omega\rangle$ exemplifies $\D$ does not have the
$M$-property. 
\end{quotation}
Clearly $T^\varepsilon_n$ is a tree of the right type, in particular
\[\mu(\lim(T^\varepsilon_n))=\inf\{|T^\varepsilon_n \cap {}^\ell 2|/
2^\ell:\ell < \omega\} \ge \varepsilon.\]
So assume $A\in\D$ and we are going to prove that $\bigcap\limits_{n \in A}
\lim(T^\varepsilon_n)$ is empty. We know that $\{\eta^\varepsilon_n:n\in A\}$
is somewhere dense, so there is $\ell^*<\omega$ and $t^*\in S^\varepsilon_{
\ell^*}$ such that:
\[\ell^*<\ell<\omega\ \&\ t^*\prec t\in S^\varepsilon_\ell\quad
\Rightarrow\quad (\exists n\in A)(t\vartriangleleft\eta^\varepsilon_n).\]
Now ${\frac{|t^* \cap {}^{\ell^*}2|}{2^{\ell^*}}}$ is $> \varepsilon$ (so
$S^\varepsilon_\ell$ was defined).  So we choose $\ell > \ell^*$, such
that:
\begin{quotation}
\noindent if $\nu\in {}^\ell 2$, $\nu\restriction\ell^*\in t^*$ 

\noindent then $t'_\nu=\{\rho\in {}^\ell 2:\rho\restriction\ell^*\in t^*\mbox{
and }\rho\ne\nu\}\in S^\varepsilon_\ell$, 
\end{quotation}
hence there is $n=n_\nu\in A$ such that $t'_\nu$ appears in
$\eta^\varepsilon_n$. Now clearly 
\[\begin{array}{ll}
\bigcap\limits_{n \in A} \lim(T^\varepsilon_n) &\supseteq \bigcap\limits_{
\scriptstyle\nu\in {}^\ell 2\atop\scriptstyle\nu\restriction\ell^*\in t^*}
\lim(T^\varepsilon_{n_\nu})\\
  &\supseteq\{\eta\in\fs:\eta\restriction\ell\in\bigcap\{t'_\nu:\nu\in {}^\rho
2,\nu\restriction\ell\in t^*\}\}=\emptyset, 
  \end{array}\]
finishing the proof.  \QED

\begin{conclusion}
\label{2.7}
In the universe $\V^{\bP_{\omega_2}}$ from \ref{2.4}, there is no
(non-principal) ultrafilter (on $\omega$) with property $M$. \QED
\end{conclusion}

\begin{concrem}
{\em 
We may consider some variants of $\qdwa$.

In definition \ref{1.1} we have $\dom(H^p)$ is as in \ref{1.1}(1) but:
$H^p\rest B^p_1$ gives constants (not functions) and for $x^m_i \in
B^p_3\setminus B^p_1$, letting $n=\min(m/E^p)$ the function $H^p(x^m_i)$
depends just on $\{x^n_j : j\leq i\}$. Moreover, it is such that changing the
value of $x^n_i$ changes the value, so $H^p(x^m_i) =x^n_i \times f^p_{x^m_i}
(x^n_0,\ldots, x^n_{i-1})$. Call this $\qtrzy$.

A second variant is when we demand the functions $f^p_{x^m_i} (x^n_0,\ldots, 
x^n_{i-1})$ to be constant, call it $\qcztery$.  

Both have the properties proved $\qdwa$. In particular, in the end of the
proof of \ref{1.7}(5), we should change: $H^r(x^m_i)$ is defined exactly as in
the proof of \ref{1.7}(4) except that when $i<h(n^*)$, $k=\min(m/E^p)$,  $k
\notin\dom(E^q)$, $k\notin \bfu$ (so $m,k,n^*$ are $E^r$--equivalent) we let
$H^r(x^k_i) = H^q(x^m_i)\times f(x^{n^*}_i) \times x^{n^*}_i$ (the first two
are constant), so $H^r(x^m_i)$ is computed as before using this value.
}
\end{concrem}

\bibliographystyle{lit-plain}
\bibliography{lista,listb,listf,listx}
\end{document}